\documentclass[12pt,a4paper]{article}
 \usepackage{amsmath,amstext,amssymb,amscd}
 \usepackage{graphicx}
\begin{document}

\title{Choice functions on posets}

\author{Danilov V. I.\thanks{ Central Institute of Economics and Mathematics of the RAS, 47, Nakhimovskii Prospect, 117418 Moscow, Russia. e-mail: vdanilov43@mail.ru }}

\maketitle

\begin{abstract}

In the paper we study choice functions on posets satisfying the conditions of heredity and outcast. For every well-ordered sequence of elements of a poset, we define the corresponding `elementary' choice function. Every such a choice function satisfies the conditions of heredity and outcast. Inversely, every choice function satisfying the conditions of heredity and outcast can be represented as a union of several elementary choice functions. This result generalizes the Aizerman-Malishevski theorem about the structure of path-independent choice functions.


\end{abstract}

\section{Introduction}

The study of choice functions began in the theory of rational decision-making. A class of choice  functions introduced by Plott \cite{P} under the name ``path-independent'' was particularly distinguished. A comprehensive description of such choice functions was  given by Aizerman and Malishevsky in \cite{AM}. Later it turned out that these functions appear naturally  in other theories, such as the theory of nonmonotonic logic \cite{KLM} and the stable contract theory \cite{AG,FFT,KC}. In the last theory, the choice functions were used to describe preferences of agents, and choice functions proposed by Plott turned out to be the most appropriate for this situation.

In particular, in the theory of contracts, the need for generalization and transfer of the concept of such choice functions to partially ordered sets (posets) was revealed. For example, the papers  \cite{AG,BB}  considered contracts that could be concluded with some intensity between 0 and 1. To ensure the existence of stable systems of contracts, the authors had to transfer the concept of path-independent choice functions to sets equipped with partial order.  Further I call such choice functions \emph{conservative}. In \cite{FFT}, the existence and good properties of stable systems of contracts under the assumption of conservativeness of choice functions of the agents were proved.

However, some important questions remained open - whether there are many such choice functions, how to construct them, what is the structure on the set of conservative choice functions? For instance, \cite{AG}  limited themselves by two examples of conservative choice functions; \cite{FFT} does not even do that. In this paper, we answer these questions. Namely, for each sequence of elements of a poset $P$, we construct the corresponding `elementary' conservative choice function on $P$. We show that an arbitrary conservative choice function is represented as the union of several elementary choice functions. This result generalizes the Aizerman-Malishevsky theorem \cite{AM} just like its infinite variant from \cite{DKS}.

In order to make the presentation more understandable, we first consider a more simple case of finite poset $P$. Then we remove the finiteness requirement. We begin with a reminder of some concepts and statements about posets and choice functions on them.

      \section{Preliminaries}

\textbf{Posets.} A \emph{poset} is a partially ordered set, i.e. a set $P$ equipped with an order relation $\le $ (reflexive, transitive, and antisymmetric) on it. Since this relation will not change, the poset is simply denoted as $P$. A poset is called \emph{linear} (or a \emph{chain}) if any two of its elements are comparable ($x\le y$ or $y\le x$). A poset is called \emph{discrete} (or trivial), if any two elements are incomparable. A more general class of posets covering two previous ones is distinguished by the transitivity condition of the comparability relation. Structurally, such a poset is a direct sum of chains. Exactly such posets were used in \cite{AG}.

A subset $I$ of $P$ is called an (order) \emph{ideal} (or a lower set, or a minor set), if with each element it contains every smaller one. For example, the principal ideal $I(x)=\{y\in P, y\le x\}$. An ideal generated by a subset $A$ in $P$ is denoted as $I(A)$; $I(A)=\cup_{a\in A} I(a)$. The dual is the concept of a \emph{filter}; a filter with each element contains larger ones. A filter generated by a subset $A$ is denoted as $F(A)$, so that $F(A)=\{x\in P, x\ge a \text{ for some } a\in A\}$. The set of all ideals is denoted by $\mathcal I(P)$; it is a complete distributive  sublattice of the Boolean lattice $2^P$ of all subsets in $P$. \medskip

\textbf{Choice functions.} In a classical situation (when a poset $P$ is trivial that is simply a non-structurated set), a choice function is a mapping $f:2^P \to 2^P$ such that $f(X)\subseteq X$ for any $X\subseteq P$. In decision theory, choice functions are used to describe behavior; having access to a variety of alternatives $X$ the decision-maker selects a subset $f(X)$. A rational decision-maker chooses the best alternatives in some sense. The rationality conditions of the corresponding choice function have been intensively studied in the choice theory, see, for example, \cite{AM}. Two conditions turned out to be the most popular. These are heredity and outcast conditions.

The \emph{heredity} property (known also as substitutability or persistency): if $A,B\subseteq P$  and $A\subseteq B$ then $f(B)\cap A\subseteq f(A)$.\medskip

In other words, if an element from a smaller set $A$ is chosen in a larger set $B$ then it should be chosen in the smaller one. \medskip

The \emph{outcast} property (known also as consistency or as Irrelevance of Rejected Alternatives): if $f(B)\subseteq A\subseteq B$ then $f (A)=f (B)$.\medskip

In words -- removing `bad' (not selected) elements does not affect the choice.\medskip

All these notions are transferred to posets without any changes. A \emph{choice function} (CF) on a poset $P$ is a mapping $f:\mathcal I (P) \to \mathcal I (P)$, such that $f(X)\subseteq X$ for any ideal $X$ in $P$. The conditions of heredity and outcast are formulated as before. \medskip

\textbf{Definition.} A CF is called \emph{conservative} if it has the heredity and outcast properties. \medskip

One can formulate the conservativeness by a single condition (\cite{KKV}): if $f(A)\subseteq B$ then $f(B)\cap A\subseteq f(A)$. We find more convenient to use heredity and outcast separately. Here are some simple properties of conservative CFs.\medskip

1. $f(f(X))=f(X)$ for any ideal $X$.\medskip

Apply the outcast to the chain of inclusions $f(X) \subseteq f(X)\subseteq X$.\medskip

2. Any conservative CF $f$ has the so-called \emph{path independence} property (or satisfies the Plott equality): for any ideals $X$ and $Y$
                                                  $$
                               f(X\cup Y)=f(f(X)\cup f(Y)).
                                                    $$

Indeed, $f(X)\subseteq X\cup Y$; from the heredity we obtain $f(X\cup Y)\cap X\subseteq f(X)$. Similarly, $f(X\cup Y)\cap Y\subseteq f (Y)$. Hence $f(X\cup Y)= f(X\cup Y)\cap (X\cup Y)\subseteq f(X)\cup f (Y)\subseteq X\cup Y$. Using the outcast, we get $f(X\cup Y)=f(f(X)\cup      f(Y))$.

In fact, this equality is true not only for two ideals, but for an arbitrary number of them. The proof is the same. Note also that the path independence implies the outcast property, but not the heredity.\medskip

3. The union of conservative CFs is a conservative CF.\medskip

Indeed, let $f=\cup _if_i$, that is $f(X)=\cup_i f_i(X)$ for any ideal $X$. We have to check the heredity and outcast properties for CF $f$.

Heredity. Let $Y\subseteq X$. Then for any $i$ we have inclusions $f_i(X)\cap Y\subseteq f_i(Y)$.       And it remains to take the union and use the distributivity.

Outcast. Let $f (X)\subseteq Y\subseteq X$. Since       $f_i (X)\subseteq f (X)$, we have the chain  $f_i (X)\subseteq Y\subseteq X$ for any $i$. The outcast for $f_i$ implies the equality  $f_i (X)=f_i (Y)$ for any $i$, from where $f(X)=f(Y)$.\medskip

This property means that having several conservative CFs $f_i$ we can form new CF $\cup_i f_i$ which also is conservative. We come to the task of finding an enough big stock of `simple' conservative CFs sufficient for construction of any conservative CF. \medskip

\textbf{Example} (a `constant' CF). Let $I$ be some ideal in $P$. For an arbitrary ideal $X$, we put $f^I (X)=I\cap X$. It is easy to see that $f^I$ is a conservative CF on $P$.\medskip

If a poset $P$ is linear, then the previous construction gives all conservative CFs. In fact, set $I=f (P)$. Since an arbitrary ideal $X$ lies in $P$, we have from heredity the inclusion $I\cap X\subseteq f (X)$. Due to the linearity of $P$, either $X\subseteq I$ or $I\subseteq X$. In the first case      $I\cap X=X$, $X\subseteq f(X)$,  from where the equality $f(X)=X=I\cap X$ comes. In the second case, $I\cap X=I$. We have the chain $I=f (P)\subseteq f(X)\subseteq P$, and the outcast gives  $f(X)=f(P)=I=I\cap X$.

However, in the case of non-linear posets, there are other conservative CFs. Two interesting examples are given in \cite{AG}. Since the union of constant CFs is constant again, we need more flexible construction of `simple' conservative CFs. To make the ideas more transparent, we will temporarily assume that poset $P$ is finite. In Section 6 we consider the general case.

      \section{Elementary CFs}

Let $A=(a_1,...,a_k)$ ($k\ge 0$) be a sequence of elements of a poset $P$ . We use this sequence to build the following CF $f_A$. Suppose that $X$ is an arbitrary ideal in $P$. We have to define $f_A(X)$. To do this denote by $i=i(A, X)$ the first number $i$, such that $a_i\in X$. In other words, $a_1,..., a_{i-1}$ does not belong to $X$, but $a_i$ does. If none $a_i$ belongs to $X$, we set $i=k$. If $k=0$ (that is if $A$  is the empty sequence) we set $i=0$.

Now, we put $f_A (X)$ to be equal to the intersection of $X$ with the ideal $I(a_1,...,a_i)$ generated by $a_1,..., a_i$, that is
                   $$
                           f_A(X)=X\cap I(a_1,...,a_i).
                     $$
In particular, $f_\emptyset(X)=\emptyset$ for any $X$. We say that the function $f_A$ is the \emph{elementary CF}, \emph{associated with the sequence} $A=(a_1,...,a_k)$.

It is easy to understand that constructing such CF, one can limit oneself to non-repeating sequences, that is, assume that all $a_i$ are different. We will do so.

The intuitive meaning of the elementary CF $f_A$ associated with the sequence $A=(a_1,...,a_k)$ is the follows. We understand the sequence $(a_1,...,a_k)$ as a hierarchy of goals of our decision-maker, and the importance of goals $a_i$ decreases with the growth of $i$. The decision-maker tries to choose the most important goal $a_1$ first of all. If it is available, that is, if it lies in the ideal $X$, he selects it (along with all smaller elements) and settles down (that is, completes the choice). If the goal $a_1$ is not available, it includes in the choice all elements of $X$ which are less than $a_1$, and proceeds to achieve the next goal $a_2$. And so on.\medskip

      \textbf{Proposition 1.} \emph{$f_A$ is a conservative CF.}\medskip

\emph{Proof.}  If the sequence $A$ is empty, then $f_A=\emptyset $ and there is nothing to check.

Let us check the heredity. Suppose $Y\subseteq X$, and $y$ belongs to both $Y$ and $f_A (X)$. We need to show that $y\in f_A (Y)$. Let $a_i$ be the first term of the sequence $A=(a_1,..., a_k)$ that falls into $X$. The elements $a_1,..., a_{i-1}$ do not belong to $X$, and therefore do not belong to $Y$. By the construction $f_A$, $y\le $ of some of $a_1,..., a_i$. But then $y$  less than the same $a_j$, and therefore belongs to $f_A (Y)$.

Let us check the outcast. Let $f_A (X)\subseteq Y\subseteq X$. It is enough to check that  $f_A(Y)\subseteq f_A(X)$. Let $a_i$ be the first element of the sequence which lies in $X$. Then $a_1,..., a_{i-1}$ do not fall in $X$, and more so in $Y$. As for $a_i$, it belongs to $f_A (X)$ and therefore  belongs to $Y$. If $y\in f_A (Y)$, then $y$ lies under some of $a_1,..., a_i$. And since $y\in X$, then it belongs to $f_A (X)$. $\Box$\medskip

\textbf{Definition.}  A sequence $A=(a_1,..., a_k)$ is \emph{compatible with} a CF $f$ if, for any $i$ from 1 to $k$, $a_i\in f(P-F(a_1,..., a_{i-1}))$. (For $i=1$, this means $a_1\in f(P)$. Recall that $F(a_1,...,a_j)$ is the filter generated by $a_1,...,a_j$.) \medskip

\textbf{Proposition 2.} \emph{If a sequence $A$ is compatible with a hereditary CF $f$ then} $f_A\subseteq f$.\medskip

\emph{Proof.} Let $X$ be an arbitrary ideal; we need to show that $f_A(X)\subseteq f(X)$. Recall how $f_A(X)$ was constructed. We find the first member $a_i$ of the sequence such that $a_i\in X$; then $f_A(X)=X\cap I(a_1,...,a_i)$. Suppose that $x$ is an element of $f_A(X)$ and $x\le a_j$ for some  $j\le i$. Since $a_j\in f(P-F(a_1,...,a_{j-1}))$ and  $f(P-F(a_1,...,a_{j-1}))$ is an ideal, we obtain that $x\in f(P-F(a_1,...,a_{j-1}))$. Because $x\in X\subset P-F(a_1,...,a_{j-1})$, the heredity of $f$ implies that $x\in f(X)$.  $\Box$\medskip

      \section{The main theorem -- a finite case}

\textbf{Theorem.} \emph{Let $f$ be a conservative CF on a finite poset $P$. Then it is the union of some elementary CFs.}\medskip

Actually we shall prove more precise assertion: \emph{any conservative CF $f$ is the union of elementary CFs associated with AC-sequences compatible with $f$.} Here a sequence $(a_1,..., a_k)$ is called \emph{antichain sequence} (AC-sequence) if all $a_i$ are incomparable in the poset $P$. This assertion follows from Proposition 2 and the following Proposition 3. \medskip

\textbf{Proposition 3.} \emph{Let $f$ be a CF on finite poset $P$ satisfying the outcast condition. Suppose that $x\in f(X)$ for some ideal $X$. Then there exists an AC-sequence $A$ compatible with $f$ such that $x\in f_A(X)$.} \medskip

\emph{Proof.} One can suppose that $f$ is non-empty CF.  We construct such a sequence $A$ step by step.

Let us discuss the first step of the construction. If $x\in f(P)$, we put $a_1=x$ and terminate the construction of the sequence. The definition of $f_A$ shows that $x\in f_A(X)$.

So we can assume that $x$ does not belong to $f(P)$. This is possible only if $f(P)$ is not contained in $X$. Indeed, otherwise $f(P)\subseteq X\subseteq P$ and from the outcast we get $f(P)=f(X)$ and $x$ belongs to $f(P)$, contrary to the assumption. Hence, the set $f(P)-X$ is non-empty. We take $a_1$ to be a minimal element in the set $f(P)-X$, put $B_1=P-F(a_1)$, and go to the second step. Note that $x\notin I(a_1)$.

At the $k$-th step, we have:

a) an AC-sequence $(a_1,..., a_k)$,

b) $x$ does not belong to the ideal $I(a_1,..., a_k)$,

c) $X\subseteq B_k:=P-F(a_1,..., a_k)$ (where $B_0=P$),

d) for any $i$ from 1 to $k$, $a_i$ is a minimal element of the set $f(B_{i-1})-X$.

In particular, d) implies that the sequence $(a_1,...,a_k)$ is compatible with $f$. As above, we consider two cases:  when $x$ belongs to $f(B_k)$ and when it doesn't.

If $x\in f(B_k)$, we set $a_{k+1}=x$ and terminate the construction of the  sequence $A$. It is clear that $x\in f_A (X)$. Due to b) and c), $a_{k+1}$ is not comparable with $a_1,...,a_k$.

Now suppose that $x\notin f(B_k)$. If $f(B_k)\subseteq X$, then from c) and the outcast we have that $f(B_k)=f (X)$ and contains $x$; a contradiction. Hence $f(B_k)$ is not contained in $X$. Put $a_{k+1}$ to be a minimal element of the non-empty set $f(B_k)-X$. We assert that  \emph{the extended sequence $(a_1,..., a_k, a_{k+1})$ also satisfies the properties a)-d). }

Let us prove a). Since $a_{k+1}\in B_k$, $a_{k+1}\notin F(a_1,...,a_k)$. Therefore we have to show that $a_{k+1}\notin I(a_1,...,a_k)$. Suppose that $a_{k+1}\le a_i$ for some $i\le k$. Since $a_i\in f(B_{i-1})$, $a_{k+1}$ belongs to $f(B_{i-1})$ as well and does not belong to $X$. Since $a_i$ is a minimal element of $f(B_{i-1})-X$ (due to d)), we obtain that $a_{k+1}=a_i$. But this is contrary to the fact that $a_{k+1}\in B_k$ and $B_k$ (see c)) does not contain $a_i$.

b) We have to show that $x$ does not lie under $a_{k+1}$. But if $x\le a_{k+1}$, then $x$ belongs to $f(B_k)$ due to the ideality of $f(B_k)$, which contradicts the assumption that $x\notin f(B_k)$.

Check c), that is $X\subseteq B_{k+1}$, or that $X$ does not intersect with the filter $F(a_{k+1})$. If there is an $y\in X$ such that $y\ge a_{k+1}$, then by ideality of $X$ the element $a_{k+1}$ also belongs to $X$, which contradicts to the choice $a_{k+1}$ outside $X$.

Finally, d) follows from the previous d) and the choice of $a_{k+1}$. \medskip

Since the poset $P$ is finite, sooner or later the process ends, and we get an AC-sequence $A$ such that $x\in f_A (X)$. $\Box$  \medskip

\section{Simplicity of elementary CFs}

We have shown that any conservative CF is represented as the union of several elementary CFs associated with AC-sequences. Now we will show that these elementary blocks are 'simple' in the sense that they no longer decompose into a union of other conservative CFs. In other words, that they are join-irreducible.\medskip

\textbf{Lemma 1.} \emph{Let $f$ and $g$ be elementary CFs associated with AC-sequences $A=(a_1,...,
a_n)$ and $B=(b_1,..., b_k)$ correspondingly, with $g\subseteq f$. Suppose, that $a_1=b_1,..., a_{i-1}=b_{i-1}$ but $a_i$ is different from $b_i$. Then $b_i<a_i$ (either $b_i$ is missing).} \medskip

Proof. We will assume that $b_i$ is actually present, and consider the ideal $X=I(a_i, b_i)$. We
state that $b_i\in g(X)$. If this is not the case, then $b_i$ is not the first member of the sequence $B$, belonging to $X$; there is a smaller number $j<i$, such that $b_j\in X$. Since it is not it is true that $b_j\le b_i$ (due to the incomparability of members of $B$), then $b_j\le      a_i$. But $b_j=a_j$, and we again get a contradiction with the incomparability of members of $A$.

Similarly $a_i$ is the first member of the sequence $A$ that falls into $X$. So
$f (X)=I(a_1,..., a_i)$. Therefore, $b_i\in I(a_1,..., a_i)$. $b_i$ cannot be less than
$a_1=b_1,..., a_{i-1}=b_{i-1}$, because this would contradict the incomparability of members of
the sequence $B$. So $b_i\le a_i$. $\Box$\medskip

\textbf{Proposition 4.} \emph{Let $f=f_A$ be an elementary CF associated with AC-sequence $A=(a_1,..., a_n)$. If $f=g\cup h$, where $g$ and $h$ are conservative CFs, then $g$ or $h$ is equal to $f$.}\medskip

Proof. Using the theorem, we decompose $g$ and $h$ into elementary CFs. As a result, we get the decomposition $f=g^1\cup ...\cup g^l$, where each CF $g^j=f_{B^j}$ is an elementary CF associated with AC-sequence $B^j=(b_1^j,..., b_{k_j}^j)$. We want to show, that at least one of $B^j$ is equal to $A$.

Consider one of these sequences $B=(b_1,..., b_k)$. For some time, it may coincide with the sequence $A$. We denote by $d(B)$ the first index $i$ for which $b_i$ is different from $a_i$ (for example,
$b_i$ is simply missing).

Assuming that all $B^j$ are different from $A$, we get that all $d(B^j)\le n$. Denote by $d$ the maximum of the numbers $d (B^1),..., d(B^l)$. We call the index $j$ a `leader' if $d (B^j)=d$.

Let us form the ideal $X=I^0(a_1)\cup ...\cup I^0(a_{d-1})\cup I (a_d)$, where $I^0 (a)=\{x\in P, x<a\}=I(a)-\{a\}$. It is clear that $a_d$ is the first member of the sequence $A$ that falls into $X$. Because $a_i$ (with $i<d$) does not belong to $I^0(a_i)$, and also does not belong to other ideals $I^0(a_j)$ due to the incomparability with the rest $a_j$. Therefore, $f (X)=X$ and, in particular, $a_d\in f(X)$.

           Let us now consider $g^j (X)$, where $j=1,..., l$; the corresponding
      the sequence $B^j$ is simply denoted as $B=(b_1,..., b_k)$. Let
      $i$ is the first number for which $b_i$ falls into X. It can't be
      a number less than $d (B)$, because for such a number $j$ (less than $i$ and so
      more than the smaller $d$) $b_j=a_j$ and does not belong to $X$ (see the lemma above). On the
other hand, $b_{d (B)}$ is less than $a_{d (B)}$ (see Lemma), and therefore
      belongs to $X$. If $j$ is not a leader, then

      $$g^j(X)=I^0(a_1)\cup ...\cup I^0(a_{d(B^j)})$$

\noindent (the last term is missing if $b_{d(B^j})$is missing). If $j$ is a leader, then

      $$g^j(X)=I^0(a_1)\cup ...\cup I^0(a_{d-1})\cup I(b^j_d)$$.

Since $b^j_d<a_d$ (by Lemma), we see that $a_d$ does not belong to any ideal $g^j (X)$. In contradiction with $a_d\in f(X)=\cup _j g^j(X)$. $\Box$

\section{The main theorem -- a general case}

We assumed above that the poset $P$ is finite. Now we will consider the general case. Everything is done the same way, only finite sequences need to be replaced by infinite ones, indexed by elements of well ordered sets (or ordinal numbers). Recall that a linearly ordered set $(I,\prec )$ is said to be \emph{well ordered} if any non-empty subset in it has a minimal element. To be not confused, ideals in the set of indices $I$ will be called \emph{initial segments}. Note that a chain $(I,\prec)$ is well ordered if and only if any initial segment of $I$, other than the entire $I$, has the form $[\prec i]=\{j\in I, j\prec i\}$ for some $i\in I$.

A \emph{well} \emph{sequence} in $P$ is a sequence $A=(a_i, i\in I)$ elements of the poset $P$, which set of indices $I$ is well ordered. With such a sequence, one can associate the CF $f_A$ on the poset $P$, which is actually defined as before. Namely, let $X$ be an ideal in $P$; denote $i=i(X)$ the first member of the sequence which belong to $X$. In other words, $a_i\in X$, and $a_j\notin X$ for $j\prec i$. Then $f_A (X)$ is equal to the intersection of $X$ with the ideal in $P$ generated by all $a_j$, $j\le i$. (If there is no such number $i$, then the ideal is generated by all $a_j$, $j\in I$.)

It is easy to understand that we can confine ourselves by non-repeating sequences, when all $a_i$ are different. In this case, one can consider $A$ as a subset of $P$ equipped with a well order $\prec$.

As in Proposition 1, it is checked that the CF $f_A$ is conservative.

The definition of compatibility of a sequence $(a_i, i\in I)$ with CF $f$ remains the same: for any $i\in I$, $a_i\in f(B_i)$, where $B_i=P - \cup _{j<i} F(a_j)$. And the same reasonings as in Proposition 2 show that $f_A\subseteq f$ if $A$ is compatible with a hereditary CF $f$.

The main theorem says now that \emph{any conservative CF $f$ is the union of some elementary CFs compatible with $f$.}

The most subtle part is a construction of well sequences compatible with CF $f$. The following reasoning resembles Zermelo's proof of that any set have a well order.

Namely, let $\mathcal F $ denote the set of filters $F$ of $P$, for which $f(P-F)\ne  \emptyset $. We fix some `selector' $p:\mathcal F  \to P$, $p(F)\in f(P-F)$. \medskip

\textbf{Definition.} A \emph{gallery} is a subset $U$ in $P$, equipped with a linear order $\prec_U$, which have the following property:\medskip

($*$) if $V$ is an initial segment in $U$ other than $U$ then $V=[\prec _U x]$, where $x=p(F(V))$.\medskip

Recall that $F(V)$ denotes the filter of $P$ generated by the set $V$. The gallery $U$ is called \emph{through} one if $F(U)$ does not belong to $\mathcal F $, that is, if $f(P-F(U))$ is empty.

By virtue of ($*$), it is obvious that the linear order $\prec_U$ of any gallery $U$ is a well order. Therefore, a gallery can be considered as a well sequence that is obviously compatible  with CF $f$. It remains to show that there are quite a lot of galleries. More precisely, we show that, \emph{for any $x\in f(X)$, there is a through gallery $U$ such that $x\in f_U(X)$.}

But first we need to say about the main property of galleries.

Say that a gallery $U'$ \emph{continues} a gallery $U$, if $U\subseteq U'$ and $U$ (as an ordered set) coincides with some initial segment $U'$. This define an order relation on the set of all galleries. Basic property of galleries is that the order is linear, i.e. any two galleries are comparable, one of them continues the other. \medskip

\textbf{Lemma 2.} \emph{Any two galleries are comparable.}\medskip

Indeed, let $U$ and $V$ be two galleries, and $W$ be the largest initial segment in both galleries, i.e. $W=\{x\in U\cap V, \ [\prec _U x]=[\prec _V x]$, and the restrictions $\prec _U$ and $\preceq _V$ on this set are the same\}.

We claim, that $W$ is equal to $U$ or $V$; this is exactly what should be proved.

Let us assume that this is not the case, and that $W$ is different from $U$ and $V$. Since $W$ is an initial segment in $U$, then, by  the property ($*$), $W$ has the form $[\prec _U x]$, where $x=p(F(W))$. Similarly, $W$ has  the form $[\prec _V y]$, where $y=p(F(W))$. So $x=y$. But then $x=y$ belongs to both $U$ and $V$, hence belongs to $W$, despite the fact that $W=[\prec _U x]$. A contradiction. $\Box$\medskip

\textbf{Corollary.}  \emph{For any selector $p$ there exists a unique through gallery.}\medskip

Proof. Take the union of all galleries. $\Box$ \medskip

\textbf{Proposition 5.} \emph{Let $X$ be an ideal of $P$, and $x\in f (X)$. Then there is a through gallery $U$ such that $x\in f_U (X)$.}\medskip

Proof. To prove this, we take a special selector $p$. Namely, assume that a filter $F$ does not intersect with the ideal $X$, that is $X\subseteq P-F$. Then there are two possible situations.
The first one is when $f (P-F)$ is not contained in $X$; in this case, we choose $p (F)$ outside of X. The second one is $f (P-F)\subseteq X$; since $X\subseteq P-F$  then from the Outcast $f (P-F)=f (X)$ and contains $x$; in this case, we choose $p(F)=x$.

Now let $U$ be a through gallery for the selector $p$ which exists due to Corollary. How does $f_U(X)$ look like? Let $V$ be the largest initial segment of the `sequence' $U$ that does not intersect with $X$. It can not be entire $U$, since then $F(U)$ does not intersect with $X$, $X\subseteq P-F(U)$, and from the outcast property  $f(X)$ is empty, contrary to $x\in f(X)$.

So $V=[\prec u]$ for $u=p(F(V))$ (see ($*$)). That is, $u$ is the first element of the sequence $U$ that belongs to $X$. According to the rule of $p$, this means that $u=x$. But in this case $x\in f_U (X)$, because $x$ belongs to both $X$ and the ideal $I(V\cup \{x\})$. $\Box$\medskip

All this together proves the main theorem in the case of an arbitrary poset $P$. \medskip

\

\noindent \textbf{Data Availability} All data generated  during this study are included in this article.

\end{document}